\documentclass[12pt]{article}
\usepackage{amssymb}
\usepackage{amsmath}
\usepackage[latin1]{inputenc} %text mit umlauten
\usepackage[T1]{fontenc}      %ec-fonts
\usepackage{graphicx}
\usepackage{eurosym}
\date{September 28, 2009}
\newtheorem{theorem}{Theorem}

\newtheorem{lemma}[theorem]{Lemma}
\newtheorem{proposition}[theorem]{Proposition}
\newtheorem{definition}[theorem]{Definition}
\newtheorem{remark}[theorem]{Remark}
\newtheorem{assumption}[theorem]{Assumption}

\newtheorem{example}[theorem]{Example}
\newcommand{\proof} [1]
{ {\bf Proof:} #1 \hfill\fbox{} \\}

\newcommand{\xad}{x_\alpha^\delta}

\newcommand{\tu}{\tilde{x}}
\newcommand{\xdag}{x^\dagger}
\newcommand{\yd}{y^\delta}

\newcommand{\req}[1]{(\ref{eq:#1})}

\newcommand{\R}{\mathbb{R}}

\newcommand\norm[1]{{\left\|#1\right\|}}

\newcommand\inner[2]{\left\langle #1,#2 \right\rangle}
\newcommand\domain[1]{{{\mathcal{D}}(#1)}}

\newcommand\set[1]{\left\{ #1 \right\}}

\title{\bf An extension of the variational inequality approach for nonlinear ill-posed problems \thanks{This paper is dedicated to Professor Charles `Chuck' Groetsch.}}

\author{{\sc Radu Ioan Bo\c{t}} and {\sc Bernd Hofmann}
\thanks{Department of Mathematics, Chemnitz University of Technology,  09107 Chemnitz,
{\sc Germany}. Email:$\;$\texttt{bot\,@\,mathematik.tu-chemnitz.de,\, hofmannb\,@\,mathematik.tu-chemnitz.de}\,.}}

\begin{document}

\maketitle

\medskip

\begin{abstract}
Convergence rates results for Tikhonov regularization of nonlinear ill-posed operator equations
in abstract function spaces require the handling of both smoothness conditions imposed on the
solution and structural conditions expressing the character of nonlinearity. Recently, the distinguished
role of variational inequalities holding on some level sets was outlined for obtaining convergence
rates results. When lower rates are expected such inequalities combine the smoothness properties of solution
and forward operator in a sophisticated manner. In this paper, using a Banach space setting we are going to extend the
variational inequality approach from H\"older rates to more general rates including the case of logarithmic
convergence rates.
\end{abstract}

\vspace{0.3cm}

{\parindent0em {\bf MSC2000 subject classification:}  47J06,}
65J20, 47A52, 49N45

\vspace{0.3cm}

{\parindent0em {\bf Keywords:} Ill-posed problems,}
Tikhonov regularization, convergence rates, variational inequalities, source conditions, distance function, structure of nonlinearity,
generalized Young inequality.

\section{Introduction}\label{s1}
\setcounter{equation}{0}
\setcounter{theorem}{0}

With the monograph \cite{Groe84} {\sc Charles Groetsch} presented an extremely well-readable introduction to the
theory of Tikhonov regularization of ill-posed operator equations in Hilbert spaces.
For linear ill-posed problems in that book the ingredients and conditions for obtaining convergence rates,
the role of source conditions and the phenomenon of saturation are outlined. The ill-posedness of a linear operator equation
describing an inverse problem with `smoothing' forward operator in Hilbert spaces corresponds with the fact that the Moore-Penrose
inverse of the forward operator is unbounded and only densely defined on the image space. In that sense, solving
linear ill-posed problems based on noisy data can be considered as the application of that unbounded operator to such data
elements. For further theoretic extensions we refer to the recent monograph \cite{Groe07}.
In 1989 {\sc Engl, Kunisch},
and {\sc Neubauer} published a seminal paper \cite{EKN89} on convergence rates results for the Tikhonov regularization of
nonlinear ill-posed problems in the Hilbert space setting (see also \cite[Chapter 10]{EHN96}).
After the turn of the millennium motivated by specific applications, for example in imaging, there occurred
numerous publications on the Banach space treatment of linear and nonlinear operator equations including convergence rates results
(see, e.g., \cite{Bremen08,GHS08,Neubau09,Resmerita05,ResmSche06,SchLouSch06}).
Initiated by the paper \cite{BurOsh04} of {\sc Burger} and {\sc Osher} Bregman distances were systematically exploited for evaluating the regularization error.
Because of a completely different methodology for obtaining convergence rates in (generalized) Tikhonov regularization we have to distinguish {\sl low rate results}
up to Bregman errors of order ${\cal O}(\delta)$ for the noise level $\delta>0$ and {\sl enhanced rate results} up to the Bregman saturation order
${\cal O}(\delta^{4/3})$. Recently, in the papers \cite{HKPS07,HofYam09}, moreover in \cite{GeHo09,HeinHof09},
in the thesis \cite{Poeschl08} and in the monograph \cite{Scherzetal09} by {\sc Scherzer et al.} the distinguished role of {\sl variational
inequalities} for proving low rate convergence rates of H\"older type was worked out.
This paper tries to extend the variational inequality approach to obtain more general Bregman rates of form ${\cal O}(\varphi(\delta))$ with concave
index functions $\varphi$. This includes the case of logarithmic convergence rates (see the papers \cite{Hohage00,Kalten08} by {\sc Hohage} and {\sc Kaltenbacher}).

The paper is organized as follows: In Section~\ref{s2} we present a general setting of Tikhonov type variational regularization with convex
stabilizing penalty functional  and strictly convex index function of the residual norm that characterize with positive some regularization
parameter  the sum functional to be minimized for obtaining stable approximate solutions of the nonlinear ill-posed problem under consideration in a
Banach space setting. The standing assumptions of the setting and assertions on weak convergence and level sets are also outlined in   Section~\ref{s2}.
The subsequent Section~\ref{s3} discusses structural conditions on the nonlinearity of the problem and source conditions as well as approximate
source conditions imposed on the solution. The first main result yielding an extension of the variational inequality approach from
convergence rates results of H\"older type to results for general convex index functions is formulated and proven as Theorem~\ref{theo1} in Section~\ref{s4}.
As an essential ingredient the proof applies a generalization of Young's inequality. The second main result will be given
in the concluding Section~\ref{s5} by the couple of Theorems~\ref{theo2} and \ref{theo3} that provide us with sufficient conditions for obtaining the
more general variational inequalities required in Theorem~\ref{theo1}. The canonical source condition for low rates in Banach spaces and distance functions
for measuring its violation form the basis for that conditions.

\section{Problem setting and assumptions}\label{s2}
\setcounter{equation}{0}
\setcounter{theorem}{0}

In this paper, ill-posed operator equations
\begin{equation}
\label{eq:ip:basic_equation} F(x)=y
\end{equation}
are under consideration,
where the operators $F:\domain{F} \subseteq X
\to Y$ with domain $\domain{F}$ are mapping between
real Banach spaces $X$ and $Y$, respectively.
For some noise level $\delta \ge 0$ let $\yd$ denote noisy data of the
exact right-hand side $y=y^0\in F(\domain{F})$ with
\begin{equation}
\label{eq:noise} \|\yd-y\|_Y \le \delta\,.
\end{equation}
Based on that data
we consider stable approximate solutions $\xad$ as minimizers
of the (generalized) Tikhonov type functional
\begin{equation}
\label{eq:ip:reg}  T_\alpha^{\delta}(x) :=
\psi(\|F(x)-\yd\|_Y)+\alpha \,\Omega(x)
\end{equation}
with a {\sl misfit function} $\psi: [0,\infty) \to [0,\infty)$ and a {\sl penalty functional}
$\Omega: \domain{\Omega} \subseteq X \to [0,\infty)$. The set of admissible elements for the
minimization of \req{ip:reg} is the intersection $\mathcal{D}:=\domain{F} \cap \domain{\Omega}$
of the occurring domains.

Index functions play a central role in our considerations. Originally coming from the theory of variable Hilbert scales
and expressing the function-valued index of such a scale element (see \cite{Hegland95,HofMat07}) we use this concept as follows:

\begin{definition} \label{def:index}
We call a real function $\eta: [0,\infty) \to [0,\infty)$ (and also its restriction to any segment $[0,\overline a]\;(0<\overline a<\infty)$) {\sl index function}
if it is continuous and strictly increasing with $\eta(0)=0$.
\end{definition}

Note that for index functions $\eta, \eta_1, \eta_2$ also the inverse function $\eta^{-1}$ and the antiderivative
$\Theta(s):=\int \limits_0^s \eta(t)dt$ are index functions, furthermore also all positive linear combinations
$\lambda_1\eta_1+\lambda_2\eta_2 \;(\lambda_1,\lambda_2 \ge 0,\, \lambda_1^2+\lambda_2^2>0)$ and compositions $\eta_1 \circ \eta_2$.

Throughout this paper we make the following assumptions:
\begin{assumption}
\label{ass1}
\begin{enumerate}
\item[] \item $X$ and $Y$ are Banach spaces with topological duals
$X^*$ and $Y^*$, respectively, where  $\|\cdot\|_X$,\; $\|\cdot\|_Y$ and  $\langle \cdot,\cdot  \rangle_{X^*,X}$ and  $\langle \cdot,\cdot  \rangle_{Y^*,Y}$ denote the associated norms and dual pairings.  In $X$ and $Y$ we consider in addition to the strong convergence $\to$ based on norms the weak convergence $\rightharpoonup$ based on the weak topology.
\item $F:\domain{F}\subseteq X \rightarrow Y$ is weakly-weakly sequentially continuous
and $\domain{F}$ is weakly sequentially closed, i.e.,
$$x_k\rightharpoonup x \quad \mbox{in} \quad X\;\;\mbox{with} \quad x_k \in \domain{F}  \quad \Longrightarrow \quad
         x \in \domain{F} \quad \mbox{and} \quad F(x_k) \rightharpoonup F(x)\quad \mbox{in} \quad Y.$$
\item The set $\domain{\Omega}$ is convex and the functional $\Omega$ is convex and weakly sequentially
lower semi-continuous.
\item The domain ${\cal D}:=\domain{F} \cap \domain{\Omega}$ is non-empty.
 \item For every $\alpha > 0,\;c \ge 0,$ and for the exact right-hand side $y=y^0$ of \req{ip:basic_equation}, the sets
      \begin{equation}
      \label{eq:ip:m_alpha}
      {\cal M}_\alpha(c):=\set{ x \in {\cal D} : T^0_\alpha(x) \leq c}
      \end{equation}
      are weakly sequentially pre-compact in the following sense:
      every sequence $\{x_k\}_{k=1}^\infty$ in ${\cal M}_\alpha(c)$
      has a subsequence, which is weakly convergent in $X$ to some element from $X$.
\item $\psi : [0,\infty) \to [0,\infty)$ is an index function with the property that there exist $a=a(\psi) > 0, b= b(\psi) > 0$ fulfilling
\begin{equation}\label{eq:psi}
\psi(u+v) \leq a \psi(u) + b \psi (v) \ \forall u,v  \in [0,
\infty).
\end{equation}
\end{enumerate}
\end{assumption}

One should notice that item 6 in Assumption \ref{ass1} is fulfilled in case $\psi : [0,\infty) \to [0,\infty)$ is a {\sl $p$-homogeneous} (with $p > 0$) and {\sl convex} index function.
We recall the $\psi$ is said to be $p$-homogeneous (with $p > 0$) whenever for all $x \in [0, +\infty)$ and all $t \geq 0$ it holds $\psi(tx) = t^p \psi(x)$.

Under the stated assumptions existence and stability of regularized
solutions $\xad$ can be shown in the lines of the proof of \cite[Theores 3.22 and 3.23]{Scherzetal09} (see also \cite[Section 3]{HKPS07}).

For the convex functional $\Omega$ with subdifferential $\partial
\Omega$ regularization errors in a Banach space setting are
frequently measured by means of Bregman distances
 \begin{equation*}
D_\xi(\tu,x) := \Omega(\tu)-\Omega(x)-\inner{\xi}{\tu-x}_{X^*,X}
\,,\quad \tu \in \domain{\Omega} \subseteq  X\,,
\end{equation*}
at $x\in \domain{\Omega} \subseteq   X$ and $\xi\in
\partial \Omega(x) \subseteq X^*$.
The set $${\mathcal D}_B( \Omega):=\set{x \in \domain{\Omega}:
\partial \Omega(x) \not=\emptyset}$$ is called Bregman domain.
An element $\xdag \in {\cal D}$ is called an {\it
$\Omega$-minimizing solution} to \req{ip:basic_equation} if
\begin{equation*}
\Omega(\xdag)= \min \set{\Omega(x) : F(x)=y, \;x \in {\cal D}} <
\infty\;.
\end{equation*}
Such $\Omega$-minimizing solutions exist under Assumption \ref{ass1} if \eqref{eq:ip:basic_equation} has a solution $x^\dag$ in ${\cal D}$.
This can be shown in analogy to the proof of \cite[Lemma 3.2]{Scherzetal09}.

We close this section by proving that the regularized solutions associated with data possessing a sufficiently small noise level $\delta$
belong to a level set like the one in  \eqref{eq:ip:m_alpha}, provided that the regularization parameters $\alpha=\alpha(\delta)$ are chosen
such that weak convergence to $\Omega$-minimizing solutions $x^\dag$ is enforced.

\begin{proposition} \label{prop0}
Consider an a priori choice $\alpha=\alpha(\delta) > 0$, $0 < \delta < \infty$, for the regularization parameter in \eqref{eq:ip:reg} depending
on the noise level $\delta$ such that
\begin{equation}\label{eq:condconv}
\alpha(\delta) \to 0 \ \mbox{and} \ \frac{\psi(\delta)}{\alpha(\delta)} \to 0.
\end{equation}
Provided that \eqref{eq:ip:m_alpha} has a solution $x^\dag$ in ${\cal D}$ then under Assumption \ref{ass1} every sequence
$\{x_n\}_{n=1}^\infty:=\{x_{\alpha(\delta_n)}\}^{\delta_n}\}_{n=1}^\infty$ of regularized solutions corresponding to a sequence
$\{y^{\delta_n}\}_{n=1}^\infty$ of data with $\lim_{n \to \infty} \delta_n = 0$ has a subsequence $\{x_{n_k}\}_{k=1}^\infty$,
which is weakly convergent in $X$, i.e. $x_{n_k} \rightharpoonup x^\dag$ and its limit $x^\dag$ is an $\Omega$-minimizing solution of
\eqref{eq:ip:m_alpha} with $\Omega(x^\dag) = \lim_{k \rightarrow \infty} \Omega(x_{n_k})$.

For given $\alpha_{\max} > 0$, let $x^\dag$ denote an
$\Omega$-minimizing solution of \eqref{eq:ip:m_alpha}. If we set
\begin{equation} \label{eq:rho1}
\rho=\alpha_{\max}(1+\Omega(x^\dag)),
\end{equation}
then we have $x^\dag \in {\cal M}_{\alpha_{\max}}(\rho)$ and there exists some $\delta_{\max} > 0$ such that
\begin{equation} \label{eq:belongM}
x^\delta_{\alpha(\delta)} \in {\cal M}_{\alpha_{\max}}(\rho) \ \mbox{for all} \ 0 < \delta \leq \delta_{\max}.
\end{equation}
\end{proposition}
\proof{The first part of the proposition can be proved in the same
manner as \cite[Theorem 3.26]{Scherzetal09}. Here the properties of
the index function $\psi$ play a determinant role.

We come now to the second part of the above statement and consider
an $\alpha_{\max} > 0$. Because of \eqref{eq:condconv} there exists
some $\delta_{\max} > 0$ such that $\alpha(\delta) \leq
\alpha_{\max}$ and $\frac{\psi(\delta)}{\alpha(\delta)} \leq
\min\{\frac{1}{2}, \frac{1}{2b}\}$ for all $0 < \delta \leq
\delta_{\max}$. In the following we write for simplicity $\alpha$
instead of $\alpha(\delta)$.

For all $0 < \delta \leq \delta_{\max}$, by \eqref{eq:psi}, it holds
$$T_{\alpha_{\max}}^0(\xad) = \psi(\|F(\xad)-y\|_Y)+\alpha_{\max}\Omega(\xad) \leq a\psi(\|F(\xad)-\yd\|_Y) + b \psi(\delta) +  \alpha_{\max}\Omega(\xad)$$
$$ = a[\psi(\|F(\xad)-\yd\|_Y) + \alpha \Omega(\xad)] + b \psi(\delta) + (\alpha_{\max} - a \alpha)\Omega(\xad) $$
$$\leq aT_\alpha^\delta(x^\dag) + b \psi(\delta) + (\alpha_{\max} - a \alpha)\Omega(\xad) \leq (a+b) \psi(\delta) + a\alpha\Omega(x^\dag) + (\alpha_{\max} - a \alpha)\Omega(\xad).$$
On the other hand, from $T_\alpha^\delta(\xad) \leq
T_\alpha^\delta(x^\dag)$ it yields $\Omega(\xad) \leq
\frac{\psi(\delta)}{\alpha} + \Omega(x^\dag)$. Consequently,
$$T_{\alpha_{\max}}^0(\xad) \leq (a+b) \psi(\delta) +
a\alpha\Omega(x^\dag) + \left(\frac{\alpha_{\max}}{\alpha} -
a\right) \psi(\delta) + \left(\alpha_{\max} -
a\alpha\right)\Omega(x^\dag)$$
$$= b\psi(\delta) +  \frac{\alpha_{\max}}{\alpha} \psi (\delta) + \alpha_{\max}\Omega(x^\dag) \leq  \alpha_{\max}(1+\Omega(x^\dag)) = \rho.$$}

\section{Source conditions and structural conditions of nonlinearity for the Banach space setting} \label{s3}
\setcounter{equation}{0} \setcounter{theorem}{0}

To obtain convergence rates for Tikhonov regularized solutions in
the case nonlinear ill-posed problems an appropriate interplay of
solutions smoothness, if possible expressed by source conditions for
$\xdag$, and of the structure of nonlinearity of $F$ in a
neighborhood of $\xdag$ is required. In this context, we are going
to restrict the situation a little bit more as follows:

\begin{assumption}
\label{ass2}
\begin{enumerate}
\item[] \item $F, \Omega,{\mathcal D}, X$ and $Y$ satisfy the Assumption~\ref{ass1}.
\item  Let $\xdag \in {\cal D}$  be  an $\Omega$-minimizing solution  of \req{ip:basic_equation}.
\item The operator $F$ is G\^ateaux differentiable in $\xdag$
with the G\^ateaux derivative \linebreak $F^\prime(\xdag) \in {\cal
L}(X,Y)$ (${\cal L}(X,Y)$ denotes the space of bounded linear operators from $X$ to $Y$).
\item  The functional $\Omega$ is G\^ateaux differentiable in $\xdag$ with the G\^ateaux derivative \linebreak $\xi=\Omega^\prime(\xdag) \in X^*$, i.e.,
$\xdag \in {\mathcal D}_B( \Omega)$ and the subdifferential
$\partial \Omega(\xdag)=\{\xi\}$ is a singleton.
\end{enumerate}
\end{assumption}

In the case of  Hilbert spaces $X$ and $Y$ by spectral theory one can consider
bounded linear operators $\eta(F^\prime(\xdag)^*F^\prime(\xdag)) \in {\cal L}(X,X)$ for any index function $\eta$
based on the fact that with the Hilbert space adjoint $F^\prime(\xdag)^*\in {\cal L}(Y,X)$ of $F^\prime(\xdag)\in {\cal L}(X,Y)$
the operators $F'(\xdag)^*F'(\xdag) \in {\cal L}(X,X)$ are non-negative and self-adjoint and this property carries over to
the operators $\eta(F^\prime(\xdag)^*F^\prime(\xdag))$.
For Banach spaces $X$ and $Y$, however, only the Banach space adjoint $F^\prime(\xdag)^* \in {\cal L}(Y^*,X^*)$ of $F^\prime(\xdag)$
is available, but $F^\prime(\xdag)^*F^\prime(\xdag)$ and hence $\eta(F^\prime(\xdag)^*F^\prime(\xdag))$ are not well-defined.
In contrast to the Hilbert space setting, where generalized source conditions
\begin{equation} \label{eq:Hilbert_general}
\xi = \eta(F^\prime(\xdag)^*F^\prime(\xdag)) v, \quad v \in X\,,
\end{equation}
can be exploited for arbitrary index functions $\eta$, in our Banach space only the source condition
\begin{equation} \label{eq:benchmark}
\xi = F^\prime(\xdag)^*\, w, \quad w \in Y^*\,,
\end{equation}
expressing a medium smoothness of $\xi$ has canonical character. We will consider this as an upper benchmark source condition
here accepting that only low and medium convergence rates for the regularized solutions are under consideration.
For expressing higher solution smoothness with respect to the stabilizing functional $\Omega$ duality mappings can be helpful
admitting enhanced convergence rates. For that we refer for example to the papers \cite{Hein09,Neubau09,NHHKT09},
but we note that the higher source conditions used there seem to be a little bit artificial.
Searching for low rate results in Banach spaces $X$ and $Y$ with solution smoothness limited by \req{benchmark}
our main drawback is the non-existence of generalized source conditions \req{Hilbert_general} with concave index functions $\eta$
such that $\sqrt{t}={\cal O}(\eta(t))$ as $t \to 0$. This class of index functions includes for $0<\nu \le 1/2$ the monomials
\begin{equation} \label{eq:monomphi}
\eta(t)\,=\, t^\nu \quad (t \ge 0)
\end{equation}
and  for all $\mu>0$ the family of logarithmic functions
\begin{equation} \label{eq:logarphi}
\eta(t)\,=\, \begin{cases}
\;0 & (t=0)\\
\;[\log(1/t)]^{-\mu} & (0<t \le e^{-\mu-1})
\end{cases}\,.
\end{equation}
Since {\sc Schock}'s paper \cite{Schock85} we know that convergence rates of regularized solutions can be arbitrarily slow. This corresponds
with arbitrarily weak solution smoothness. For example the very low multiple logarithmic rates for associated generalized source conditions with
index function $\eta(t)= \log \log ... \log(1/t)$ really
occur in applications of the Hilbert space theory.

One way of compensating the Banach space drawback of missing
generalized source conditions consists in applying the {\sl method
of approximate source conditions} (see \cite{DHY07,HeinHof09})
whenever $\xi$ fails to satisfy the benchmark source condition
\req{benchmark} for any $w \in Y$, but the obviously non-negative and
non-increasing distance function
\begin{equation} \label{eq:distance}
d(R):=\inf \{\|\xi-F^\prime(\xdag)^*\, w\|_{X^*}:\; w \in
Y^*,\;\|w\|_{Y^*} \le R \} \,,
\end{equation}
which is well-defined for all $R \ge 0$, fulfills the limit condition
\begin{equation} \label{eq:limit}
 \lim \limits _{R \to \infty} d(R)\,=\,0\,.
\end{equation}
We notice that on the one hand $d : [0,\infty) \rightarrow [0,\infty)$ is continuous, being convex, and on the other hand \eqref{eq:limit} is fulfilled if and only if $\xi \in
\overline{{\cal R}(F^\prime(\xdag)^*)}^{\|\cdot\|_{X^*}}$. By a
separation theorem one can prove that the latter is guaranteed provided
$F^\prime(\xdag)^{**}$ is injective (respectively, $F^\prime(\xdag)$
is injective, if $X,Y$ are reflexive Banach spaces). For $A \in
{\cal L}(X,Y)$ we denote by $A^{**} \in {\cal L}(X^{**},Y^{**})$,
defined by $\langle A^{**}x^{**}, y^* \rangle_{Y^{**},Y^*} = \langle
x^{**}, A^{*}y^*\rangle_{X^{**}, X^*}$ for $x^{**} \in X^{**}$ and
$y^* \in Y^*$, its {\sl bi-adjoint operator}.

The following lemma will be used in order to guarantee that the
distance function defined in \eqref{eq:distance} strictly
decreasing.

\begin{lemma} \label{lemma1}
Let $X,Y$ be reflexive Banach spaces and $A \in {\cal L}(X,Y)$ an
injective operator. For $\xi \in X^*$ we assume that $\xi \notin
{\cal R}(A^*)$. Then the distance function $d:[0,+\infty)
\rightarrow (0,+\infty)$, defined by
$$d(R) = \inf\{\|\xi-A^*w\|_{X^*} :\;w \in Y^*,\; \|w\|_{Y^*} \leq R\},$$
is strictly decreasing.
\end{lemma}
\proof{First let us notice that for all $R \geq 0$ there exists $\bar w
\in Y^*$, $\|\bar w\|_{Y^*} \leq R$, such that $d(R) = \|\xi-A^*\bar
w\|_{X^*}$. This is because of the fact that the dual norm function
is weak$^*$ lower semicontinuous and the unit ball in $Y^*$ is
weak$^*$ compact (Theorem of Alaoglu-Bourbaki).

Let be $R \geq 0$. Next we prove that if for $\bar w \in Y^*$, $\|\bar
w\|_{Y^*} \leq R$, it holds $d(R) = \|\xi-A^*\bar w\|_{X^*}$, then
one necessarily must have $\|\bar w\|_{Y^*} = R$. In case $R > 0$, this fact is obvious.
Suppose now that $R> 0$. Indeed, in this
case $\bar w$ is an optimal solution of the convex optimization
problems
$$\inf_{\|w\|_{Y^*}-R \leq 0} \|\xi - A^*w\|_{X^*}.$$
As the Slater constraint qualification is fulfilled (for $w'=0$ we
have $\|w'\|_{Y^*}-R < 0$), there exists a Lagrange multiplier $\bar
\lambda \geq 0$ such that (see, for instance, \cite[Theorem
2.9.2]{Zalinescu})
$$\bar \lambda (\|\bar w\|_{Y^*}-R) = 0 $$
and
$$0 \in \partial(\|\xi-A^*(\cdot)\|_{X^*} + \bar \lambda(\|\cdot\|_{Y^*} - R))(\bar w).$$
If we prove that $\bar \lambda > 0$, then the assertion follows. We
assume the contrary. This means that
$$0 \in \partial(\|\xi-A^*(\cdot)\|_{X^*})(\bar w).$$
Next we evaluate the above subdifferential. Let be $L :X^*
\rightarrow \R$, $L(w) = \|\xi + w\|_{X^*}$. Since $L$ is
continuous, by \cite[Theorem 2.8.2]{Zalinescu} we have that
$$\partial(\|\xi-A^*(\cdot)\|_{X^*})(\bar w) = \partial (L \circ (-A^*))(\bar w) = -A (\partial L (-A^* \bar w)).$$
As $A$ is injective,
\begin{equation}\label{eq:partial}
0 \in \partial L (-A^* \bar w) = \partial \|\cdot\|_{X^*}(\xi -A^*
\bar w).
\end{equation}
For the subdifferential of the norm we have the following
expressions
$$\partial \|\cdot\|_{X^*}(v) = \{u \in X : \|u\|_X \leq 1\}, \ \mbox{if} \ v = 0,$$
while
$$\partial \|\cdot\|_{X^*}(v) = \{u \in X : \|u\|_X = 1, \langle v, u \rangle_{X^*,X} = \|v\|_{X^*}\}, \ \mbox{if} \ v \neq 0.$$
By \eqref{eq:partial} it follows that only the first situation is
possible. Consequently, $\xi -A^* \bar w = 0$. But this is a
contradiction to $\xi \notin {\cal R}(A^*)$. Thus $\bar \lambda > 0$
and, so, $\|\bar w\|_{Y^*}=R$.

Let us prove now that $d$ is strictly decreasing. To this aim take
$R_1, R_2 \in [0, +\infty)$ such that $0 \leq R_1 < R_2$. It holds $d(R_1)
\geq d(R_2)$. Assume that $d(R_1) = d(R_2)$. Then there exists $w_1,
w_2 \in Y^*$, $\|w_1\|_{Y^*} = R_1$, $\|w_2\|_{Y^*} = R_2$, such
that $d(R_1) = d(R_2) = \|\xi-A^*w_1\|_{X^*} =
\|\xi-A^*w_2\|_{X^*}$. As $\|w_1\|_{Y^*} < R_2$, this leads to a
contradiction to the above considerations. Consequently, $d(R_1) >
d(R_2)$ and this concludes the proof.}

Let us mention that the decay rate of $d(R) \to 0$ as $R \to
\infty$, as assumed in \eqref{eq:limit}, expresses for the element
$\xi$ the degree of violation of \req{benchmark} and thus it can be
handled as a replacement information for the missing index function
$\eta$ from \req{Hilbert_general} in the Banach space setting.

As an adaption of the {\sl local degree of nonlinearity} introduced
for a Hilbert space setting in \cite[Definition~1]{HofSch94} to the
Banach space situation with Bregman distance we suggested in
\cite[Definition~3.2]{HofYam09} a definition, which attains here
under Assumption~\ref{ass2} the form:

\begin{definition} \label{def:degree}
Let $ 0 \le c_1,c_2 \le 1$ and $0<c_1+c_2 \le 1$. We define $F$ to be
nonlinear of degree $(c_1,c_2)$ at $\xdag$ for the Bregman distance
$D_\xi(\cdot,\xdag)$ of $\Omega$ with $\xi=\Omega^\prime(\xdag)$ if there is a constant $K>0$  such that
 \begin{equation}
      \label{eq:degree}
\norm{F(x)-F(\xdag)-F^\prime(\xdag)(x-\xdag)}_Y  \,\leq \,
K\;\norm{F(x)-F(\xdag)}_Y^{\,c_1}\; D_\xi(x,\xdag)^{\,c_2}
      \end{equation}
      for all $x \in \mathcal{M}_{\alpha_{max}}(\rho)$.
\end{definition}

In \cite{HeinHof09} it was shown that the method of approximate source conditions yields convergence rates
of Tikhonov regularized solutions $\xad$ minimizing \req{ip:reg} with misfit function $\psi(t)=t^p \;(p>1)$ whenever we have $c_1>0$ in the degree of nonlinearity and if $\xi$ fails to satisfy the benchmark source condition \req{benchmark}. The corresponding rates depend on the
distance function \req{distance}. If $c_1>0$ and the source condition \req{benchmark} holds, then we even obtain
H\"older convergence rates with H\"older exponents $\kappa=\frac{c_1}{1-c_2}$ (see \cite{HofYam09}). If the nonlinearity
of $F$ at $\xdag$ is such that $c_1>0$ cannot be satisfied, then rate results are only known for $c_2=1$ and \req{benchmark}
under the additional smallness condition $K\|w\|_{Y^*}<1$ (see \cite{Scherzetal09}). As already mentioned in \cite{Kalten08}
for the Hilbert space setting, there seem to be no rate assertions for $c_1=0, c_2=1$ if $\xi$ fails to satisfy the benchmark source condition \req{benchmark},
i.e., low rate results including results on logarithmic convergence rates are missing for $c_1=0$ in case of absence of the source condition \req{benchmark}.
However, if we cannot find a $c_1>0$ it is an interesting open question whether low rate results can be derived if the structure of nonlinearity only satisfies the weaker condition
\begin{equation}
      \label{eq:newstructure}
\norm{F^\prime(\xdag)(x-\xdag)}_Y  \,\leq \,
C\; \sigma(\norm{F(x)-F(\xdag)}_Y)
      \end{equation}
      for all $x \in \mathcal{M}_{\alpha_{max}}(\rho)$ with some constant $C>0$ and some  index function $\sigma$.
      We will attack this task in the next section using variational inequalities
as tool. In this context, let us note that the validity of \req{newstructure} for $\sigma(t)=t^{c_1}$
      and $0<c_1 \le 1$ implies with the triangle inequality that we have
$$\norm{F(x)-F(\xdag)-F^\prime(\xdag)(x-\xdag)}_Y  \leq \norm{F^\prime(\xdag)(x-\xdag)}_Y + \norm{F(x)-F(\xdag)}_Y $$ $$\le C\, \norm{F(x)-F(\xdag)}^{c_1}_Y + \norm{F(x)-F(\xdag)}_Y \le K\,\norm{F(x)-F(\xdag)}_Y^{\,c_1}$$
on the associated level sets which shows a degree $(c_1,0)$ of nonlinearity. In general we conjecture that only {\sl concave} index functions
$\sigma$ are of interest in the condition \req{newstructure}.

\section{Variational inequalities and convergence rates} \label{s4}
\setcounter{equation}{0} \setcounter{theorem}{0}

In recent publications (see
\cite{HeinHof09,HKPS07,HofYam09,Scherzetal09}) variational
inequalities of the form
\begin{equation} \label{eq:vareqkappa}
\inner{\xi}{\xdag-x}_{X^*,X}
      \leq \beta_1 D_\xi(x,\xdag) \, + \,
      \beta_2\norm{F(x)-F(\xdag)}_Y^\kappa \qquad \mbox{for all} \quad  x \in {\cal M}_{\alpha_{max}}(\rho)
\end{equation}
with two multipliers $0 \le \beta_1 <1, \;\beta_2 > 0$ and an
exponent $\kappa>0$ have been exploited for obtaining convergence
rates in Tikhonov regularization in Banach spaces, where in the
functional \req{ip:reg} to be minimized the strictly convex misfit
function $\psi(t)=t^p\;(p > 1)$ was used. We repeat in our context
the Proposition~3.3  from \cite{HofYam09}:

\begin{proposition} \label{prop1}
Set $\psi(t):=t^p\;(p>1)$ in \req{ip:reg} and assume that $F,
\Omega,{\mathcal D}, X, Y, \xdag$ and $\xi$ satisfy the
Assumption~\ref{ass2}. If there exist constants $0 \le \beta_1 <1,
\;\beta_2 >0,$ and $0<\kappa \le 1$ such that the variational
inequality \req{vareqkappa} holds with $\rho$ from \req{rho1}, then
we have the convergence rate
 \begin{equation}
      \label{eq:prop1rate}
D_\xi(x_{\alpha(\delta)}^\delta,\xdag)\,=\,{\mathcal
O}\left(\delta^\kappa\right) \quad \mbox{as}\quad \delta \to 0
\end{equation}
for an a priori parameter choice $\alpha(\delta) \asymp  \delta^{p-\kappa}$.
\end{proposition}
The proof of this proposition is based on the inequality
$T^\delta_\alpha(\xad) \le T^\delta_\alpha(\xdag)$ that holds for
all regularized solutions $\xad$ and on the variant
\begin{equation} \label{eq:youngeps}
a\,b\, \le a^{p_1}+\frac{b^{p_2}}{{p_1}^{p_2/p_1} p_2}\qquad (a,b
\ge 0,\;\;p_1,p_2 > 1 \;\; \mbox{with} \;\;
\frac{1}{p_1}+\frac{1}{p_1}=1)
\end{equation}
of Young's inequality. Note that due to Proposition~4.3 in \cite{HofYam09} the case $\kappa>1$
is not of interest, since \req{vareqkappa} with $\kappa>1$ implies the singular case $\xi =0$.

For obtaining more general low order convergence rates we change
\req{vareqkappa} as follows: We assume that there holds a
variational inequality
\begin{equation} \label{eq:vareqnew}
\inner{\xi}{\xdag-x}_{X^*,X}
      \leq \beta_1 D_\xi(x,\xdag) \, + \,
      \beta_2\,\varphi(\norm{F(x)-F(\xdag)}_Y) \qquad \mbox{for all} \quad  x \in {\cal M}_{\alpha_{max}}(\rho)
\end{equation}
with two multipliers $0 \le \beta_1 <1, \;\beta_2 > 0$ and an
index function $\varphi$.

\pagebreak

\begin{assumption}
\label{ass3}
Regarding the functions $\psi$ from \req{ip:reg} and $\varphi$ from \req{vareqnew} we make the following
assumptions:
\begin{enumerate}
\item $\psi$ and $\varphi$ are index
functions which are twice differentiable on the interior of their
domains.
\item  $\psi$ is strictly convex and $\varphi$ is concave.
\end{enumerate}
\end{assumption}

Under Assumption \ref{ass3} we can define another index function $f$
as follows:
\begin{equation}\label{eq:f}
f(0) = 0 \ \mbox{and} \ f(s) = \left[\frac{\psi'}{\varphi'} \circ
\varphi^{-1}\right](s) \ \mbox{when} \ s > 0.
\end{equation}
Let us first show that $f$ is well-defined, by proving that
$\varphi'(s) > 0$ when $s > 0$. Indeed, suppose that there exists
$\bar s > 0$ in the interior of the domain of $\varphi$ such that
$\varphi'(\bar s) = 0$. Take $t > \bar s$. By the concavity
assumption one has
$$0=\varphi'(\bar s)(t- \bar s) \geq \varphi(t)-\varphi(\bar s),$$
which contradicts the fact that $\varphi$ is strictly increasing.

By employing similar arguments, since $\psi$ is convex, whenever $s
>0$ one has that $\psi'(s) > 0$  and so $f(s) > 0$. As the
continuity of $f$ is automatically satisfied, in order to prove that
$f$ is an index function, one only needs to show  that $f$ is
strictly increasing.

Take $0 < s_1 < s_2$. Then $\varphi^{-1}(s_1) < \varphi^{-1}(s_2)$.
As $\psi$ is strictly convex, $\psi'$ is strictly increasing and so
$0 < \psi'(\varphi^{-1}(s_1)) < \psi'(\varphi^{-1}(s_2))$. On the
other hand, since $\varphi$ is concave, $\varphi'$ is
non-increasing, consequently, $\varphi'(\varphi^{-1}(s_1)) \geq
\varphi'(\varphi^{-1}(s_2)) > 0$. From here one has $f(s_1) <
f(s_2)$. Hence $f$ is an index function and so is the antiderivative
\begin{equation} \label{eq:H}
H(s):=\int \limits_0^s f(\tau) d\tau\,.
\end{equation}

From \eqref{eq:f} it follows that for $s> 0$ $\psi(s) =
\int\limits_{0}^{\varphi(s)} f(t)dt + C$. As $\psi(0) = 0$, this
yields $C=0$ and, consequently,
$$\psi(s) = H(\varphi(s)) =\int\limits_{0}^{\varphi(s)} f(t)dt\,.$$
Now aspects of the interplay between $\psi, \varphi, f$ and $H$ can be written in different manner by the equations
$$ \psi=H \circ \varphi, \quad  H=\psi \circ \varphi^{-1}$$
and
$$f(s)=[\psi \circ \varphi^{-1}](s)^\prime\;\; (s > 0),$$
where the last equation yields \req{f} by differentiation and use of the chain rule.
Further, let
\begin{equation} \label{eq:G}
G(s):=\int \limits_0^s f^{-1}(\tau) d\tau
\end{equation}
be the antiderivative  of the inverse function to $f$ and one can verify the cross-connections
\begin{equation}\label{eq:interplay}
H=G \circ f \quad \mbox{and} \quad G^{-1} \circ \psi= f \circ \varphi\,.
\end{equation}

Now we are ready to present the main convergence rate result of this paper:

\begin{theorem} \label{theo1}
Assume that $F, \Omega,{\mathcal D}, X, Y, \xdag$, $\xi$ and $\psi$
satisfy Assumption~\ref{ass2} and assume that $\psi$ and $\varphi$ satisfy
Assumption~\ref{ass3} which ensures the existence of an index function $f$
defined by \req{f}. Let there exist constants $0 \le \beta_1 <1,
\;\beta_2 >0,$ such that the variational inequality \req{vareqnew}
holds with $\rho$ from \req{rho1}. Then we have the convergence rate
of Tikhonov regularized solutions
 \begin{equation}
      \label{eq:theo1rate}
D_\xi(x_{\alpha(\delta)}^\delta,\xdag)\,=\,{\mathcal
O}(\varphi(\delta)) \quad \mbox{as}\quad \delta \to 0
\end{equation}
for an a priori parameter choice
\begin{equation} \label{eq:apriori}
\alpha(\delta) = \frac{1}{a \beta_2} f(\varphi(\delta)).
\end{equation}
\end{theorem}
\proof{For all $\alpha>0$ regularized solutions $\xad$ minimizing \req{ip:reg} have to satisfy
the inequalities  $T^\delta_\alpha(\xad) \le T^\delta_\alpha(\xdag)$. Using the definition of the Bregman distance this implies for the
noise model \req{noise} the estimate
\begin{equation}\label{eq:fund_est}
\psi(\norm{F(\xad)-\yd}_Y) + \alpha D_\xi(\xad,\xdag) \leq
\psi(\delta) + \alpha
\left(\Omega(\xdag)-\Omega(\xad)+D_\xi(\xad,\xdag) \right)\;.
\end{equation}
Moreover, from the variational inequality \req{vareqnew} we obtain
that
\begin{equation*}
\begin{aligned}
~& \Omega(\xdag)-\Omega(\xad) +  D_\xi(\xad,\xdag) = -\inner{\xi}{\xad-\xdag}_{X^*,X} \\
\leq& \; \beta_1\, D_\xi(\xad,\xdag) \, +
      \, \beta_2\,\varphi(\norm{F(\xad)-F(\xdag)}_Y).
\end{aligned}
\end{equation*}
Therefore from  \req{fund_est} it follows that
\begin{eqnarray} \label{eq:new_est}
\psi(\|F(x_\alpha^\delta) - y^\delta\|_Y) + \alpha
D_\xi(x_\alpha^\delta, x^\dag) & \leq & \psi(\delta) + \alpha\beta_1
D_\xi(x_\alpha^\delta, x^\dag) + \alpha \beta_2
\varphi(\|F(x_\alpha^\delta) - F(x^\dag)\|_Y)\nonumber\\  & & \\
& = & \psi(\delta) + \alpha\beta_1 D_\xi(x_\alpha^\delta, x^\dag) +
\frac{1}{a}(\alpha a \beta_2)\varphi(\|F(x_\alpha^\delta) -
F(x^\dag)\|_Y)\,. \nonumber
\end{eqnarray}
Using the generalization of Young's inequality (see, for instance,
\cite{los}) for the index function $f$ we obtain for sufficiently
small $\alpha
> 0$
\begin{eqnarray} \label{eq:youngnew}
(\alpha a \beta_2)\varphi(\|F(x_\alpha^\delta) - F(x^\dag)\|_Y) &
\le & \int \limits _0^{\varphi(\|F(x_\alpha^\delta) -
F(x^\dag)\|_Y)}f(t)dt +  \int \limits _0^{\alpha a \beta_2}
f^{-1}(\tau)d\tau\nonumber\\ & & \\
=H(\varphi(\|F(x_\alpha^\delta) - F(x^\dag)\|_Y)+ G(\alpha a \beta_2)& = & \psi(\|F(x_\alpha^\delta) -
F(x^\dag)\|_Y) + G(\alpha a \beta_2). \nonumber
\end{eqnarray}
From \eqref{eq:new_est} and \eqref{eq:youngnew} it follows that

\pagebreak

$$\psi(\|F(x_\alpha^\delta) - y^\delta\|_Y) + \alpha D_\xi(x_\alpha^\delta, x^\dag)$$
$$\leq \psi(\delta) + \alpha\beta_1 D_\xi(x_\alpha^\delta, x^\dag) + \frac{1}{a} \,\psi(\|F(x_\alpha^\delta) - F(x^\dag)\|_Y) + \frac{1}{a}\, G(\alpha a \beta_2)$$
$$\leq \psi(\delta) + \alpha\beta_1 D_\xi(x_\alpha^\delta, x^\dag) + \psi(\|F(x_\alpha^\delta) - y^\delta\|_Y) + \frac{b}{a} \psi(\delta) + \frac{1}{a}\, G(\alpha a \beta_2).$$
Consequently,
\begin{equation} \label{eq:equilibrate}
D_\xi(x_\alpha^\delta, x^\dag) \leq \frac{1}{(1-\beta_1)a} \frac{(a+b)\psi(\delta) + G(\alpha a \beta_2)}{\alpha},
\end{equation}
for sufficiently small $\alpha>0$.

Next we prove that for $\alpha(\delta) := \frac{1}{a \beta_2}
f(\varphi(\delta))$ it holds
\begin{equation}\label{eq:estimatephi}
\frac{(a+b)\psi(\delta) + G(\alpha(\delta) a
\beta_2)}{\alpha(\delta)} \leq (a+b+1)a\beta_2 \varphi(\delta)
\end{equation}
for sufficiently small $\delta > 0$. Indeed, \eqref{eq:estimatephi}
is equivalent to
\begin{equation}\label{eq:estimatephi2}
(a+b+1)a\beta_2 \varphi(\delta) \alpha(\delta) - (a+b)\psi(\delta) -
G(\alpha(\delta) a \beta_2) \geq 0
\end{equation}
for sufficiently small $\delta > 0$. Denote by $K(\delta) :=
(a+b+1)a\beta_2 \varphi(\delta) \alpha(\delta) - (a+b)\psi(\delta) -
G(\alpha(\delta) a \beta_2)$. One has $K(0) = 0$. We prove that
$K'(\delta) > 0$, for sufficiently small $\delta > 0$, and this will
have as consequence the fact that $K(\delta) > K(0) = 0$, for
sufficiently small $\delta > 0$. Indeed, one has for sufficiently
small $\delta > 0$
$$K'(\delta) = (a+b+1) \varphi'(\delta) f(\varphi(\delta)) + (a+b+1) \varphi(\delta) f'(\varphi(\delta)) \varphi'(\delta)$$
$$- (a+b) \varphi'(\delta) f(\varphi(\delta)) - f^{-1}(f(\varphi(\delta))) f'(\varphi(\delta)) \varphi'(\delta)$$
$$= \varphi'(\delta) f(\varphi(\delta)) + (a+b) \varphi(\delta) f'(\varphi(\delta)) \varphi'(\delta) > 0.$$
Thus \eqref{eq:estimatephi} holds and this yields the estimate
$$D_{\xi}(x_{\alpha(\delta)}^\delta,\xdag) \leq c_0 \varphi(\delta)$$
for sufficiently small $\delta>0$ and some constant $c_0>0$.}

We should note here that because of \req{interplay} $\alpha(\delta) = \frac{1}{a \beta_2}
f(\varphi(\delta))=\frac{1}{a \beta_2} G^{-1}(\psi(\delta))$  denotes an equilibration up to a constant of the two terms in the
numerator of the second fraction in \req{equilibrate}. This order equilibration corresponds with the standard approach for obtaining such convergence rates.

\begin{example}\label{example}
{\rm We conclude this section with the example situation of monomials (power functions) $\varphi(t)=t^\kappa \;(0<\kappa \le 1)$ and $\psi(t)=t^p\;(p>1)$
discussed in \cite{HofYam09} for which Proposition~\ref{prop1} was repeated above. Then our assumptions are satisfied and we have $$H(t)=t^{p/\kappa},\;\;
f(t)=\frac{p}{\kappa}\,t^{(p-\kappa)/\kappa}, \;\; G(t) \sim t^{p/(p-\kappa)}, \;\; \alpha(\delta) \sim \delta^{p-\kappa},
\;\; D_\xi(x_{\alpha(\delta)}^\delta,\xdag)\,=\,{\mathcal
O}(\delta^\kappa)\,.$$
We would like to notice that one comes to the same conclusion also in the case $0 < \kappa < p \leq 1$ discussed in \cite{GeHo09}. The reason therefore lay in the fact that $f$ remains an index function and, consequently, Theorem \ref{theo1} is still applicable, even if in this situation $\psi$ fails to be strictly concave. In fact, in order to obtain the convergence rate \eqref{eq:theo1rate} in Theorem \ref{theo1} one needs only to guarantee that the function $f$ defined as in \eqref{eq:f} is an index function which is differentiable on the interior of its domain. This happens when Assumption \ref{ass3} is satisfied, but can be the case also in other settings.}
\end{example}

\section{Variational inequalities based on canonical source conditions and approximate source conditions} \label{s5}
\setcounter{equation}{0} \setcounter{theorem}{0}

In this section we are going to formulate sufficient conditions for variational inequalities \req{vareqnew} when only some weak structural assumption of the form \req{newstructure} on the nonlinearity of $F$ with concave index function $\sigma$ is imposed.

\begin{theorem} \label{theo2}
Assume that $F, \Omega,{\mathcal D}, X, Y, \xdag, \xi$ and $\psi$
satisfy the Assumption \ref{ass2}. \linebreak Let $\xi$ satisfy the canonical
source condition \req{benchmark} and the structural condition
\req{newstructure} with some index function $\sigma$ and some
constant $C>0$ for all $x \in \mathcal{M}_{\alpha_{max}}(\rho)$.
Then a variational inequality \req{vareqnew} holds with two
multipliers $0 \le \beta_1 <1, \;\beta_2 > 0$ and with the index function $\varphi=\sigma$.
\end{theorem}

\proof{Owing to \req{benchmark} and \req{newstructure} we can estimate for all $x \in \mathcal{M}_{\alpha_{max}}(\rho)$ as
$$\inner{\xi}{\xdag-x}_{X^*,X}  = \inner{F^\prime(\xdag)^*w}{\xdag-x}_{X^*,X} = \inner{w}{F^\prime(\xdag)(\xdag-x)}_{Y^*,Y}$$
$$ \le \norm{w}_{Y^*} \norm{F^\prime(\xdag)(x-\xdag)}_Y \le C\,\norm{w}_{Y^*} \sigma(\norm{F(x)-F(\xdag)}_Y)\,. $$
This, however, yields the variational inequality \req{vareqnew} with $\beta_1=0<1,\;\beta_2=C\,\norm{w}_{Y^*}$ and with $\varphi=\sigma$,
where $\sigma$ is the index function
from \req{newstructure}. This proves the theorem.
}

\begin{theorem} \label{theo3}
Assume that $X,Y$ are reflexive Banach spaces,
$F, \Omega,{\mathcal D}, X, Y, \xdag, \xi$ and $\psi$ satisfy the Assumption~\ref{ass2},
and $F^\prime(\xdag)$ is an injective operator.
Let $\xi \notin {\cal
R}(F^\prime(\xdag)^*)$. Moreover, assume that the structural
condition \req{newstructure} is fulfilled with some index function
$\sigma $ and some constant $C>0$ for all $x \in
\mathcal{M}_{\alpha_{max}}(\rho)$ and that the Bregman distance is
locally $q$-coercive with $2 \le q <\infty$, i.e.~there is some
constant $c_q>0$ such that
\begin{equation} \label{eq:coercive}
 D_{\xi}(x,\xdag) \ge c_q\,\norm{x-\xdag}^q_X
\end{equation}
holds for all $x \in \mathcal{M}_{\alpha_{max}}(\rho)$.
Then a variational inequality \req{vareqnew}
 holds for all \linebreak $x \in {\cal M}_{\alpha_{max}}(\rho)$ with two multipliers $0 \le \beta_1 <1, \;\beta_2 > 0$ and
with the index function $\varphi(0)=0,\;\varphi(t)=\left[d \left(\Psi^{-1}(\sigma(t))
\right)\right]^{q^*}\;(t>0)$, where $\frac{1}{q}+\frac{1}{q^*}=1$ and $\Psi:(0,\infty) \rightarrow (0,\infty)$,
$\Psi(R):=\frac{d(R)^{q^*}}{R}$.
\end{theorem}

\proof{Instead of \req{benchmark} we have here for all $R >0$ the
equations $\xi=F^\prime(\xdag)^*w_R+r_R$ with $\norm{w_R}_{Y^*} \le
R$ and $\norm{r_R}_{X^*}=d(R)$. Then we can estimate for all $x \in
\mathcal{M}_{\alpha_{max}}(\rho)$ by using \req{newstructure} as
$$\inner{\xi}{\xdag-x}_{X^*,X}  = \inner{F^\prime(\xdag)^*w_R+r_R}{\xdag-x}_{X^*,X}\leq \inner{w_R}{F^\prime(\xdag)(\xdag-x)}_{Y^*,Y}+\inner{r_R}{\xdag-x}_{X^*,X}$$
$$ \le R\,  \norm{F^\prime(\xdag)(x-\xdag)}_Y + d(R)\,\norm{x-\xdag}  \le R\,C\, \sigma(\norm{F(x)-F(\xdag)}_Y)+d(R)\,\norm{x-\xdag} \,. $$
Now for $q$ and $q ^*$ adjoint exponents with $1/q+1/q^*=1$ the inequality
$$\inner{\xi}{\xdag-x}_{X^*,X} \le R\,C\, \sigma(\norm{F(x)-F(\xdag)}_Y)+c_q^{-1/q}\,d(R)\, D_{\xi}(x,\xdag)^{1/q}$$
obtained from \req{coercive} can be further handled by using Young's inequality in the standard form
$$a\,b\, \le \frac{a^{p_1}}{p_1}+\frac{b^{p_2}}{p_2}\qquad (a,b
\ge 0,\;\;p_1,p_2 > 1 \;\; \mbox{with} \;\;
\frac{1}{p_1}+\frac{1}{p_1}=1)$$
when setting $a:=D_{\xi}(x,\xdag),\;b:=c_q^{-1/q}\,d(R),\;p_1:=q,\;p_2:=q^*.$
In that way we derive for all $R>0$
$$\inner{\xi}{\xdag-x}_{X^*,X} \le R\,C\, \sigma(\norm{F(x)-F(\xdag)}_Y)+\frac{1}{q}\,D_{\xi}(x,\xdag)+ \frac{c_q^{-q^*/q}}{q^*}\,d(R)^{q^*}\,.$$
The continuity of $d$ carries over to the auxiliary function $\Psi:(0,\infty) \rightarrow (0,\infty)$, $\Psi(R)=\frac{d(R)^{q^*}}{R}$, which is continuous and strictly
decreasing, and which fulfills $\lim_{R \to 0} \Psi(R)=\infty$ and $\lim_{R \to \infty} \Psi(R)=0$.  Its inverse $\Psi^{-1} :(0,\infty) \rightarrow (0,\infty)$ is also
continuous and strictly decreasing and for all $t>0$ the equation $\Psi(R)=\sigma(t)$
has a uniquely determined solution $R>0$.  Note that for rates results only sufficiently small $t>0$ are of interest.
Setting $R:=\Psi^{-1}\left(\sigma(\norm{F(x)-F(\xdag)}_Y)\right)$ we get some constant $\hat C>0$ such that the variational inequality
$$\inner{\xi}{\xdag-x}_{X^*,X} \le \frac{1}{q}\,D_{\xi}(x,\xdag)+ \hat C\, \left[d \left(\Psi^{-1}(\sigma(\norm{F(x)-F(\xdag)}_Y))  \right)\right]^{q^*} $$
holds for all $x \in \mathcal{M}_{\alpha_{max}}(\rho)$. Now the function defined by $\zeta(s):=d \circ \Psi^{-1}\circ \sigma(s)$ when $s > 0$ with extension $\zeta(0):=0$
is an index function. Namely, $\zeta$ is continuous on $(0,\infty)$, since $d$ is continuous. Moreover, the limit $\lim_{R \to \infty} d(R)=0$ implies $\lim_{t \to 0} \zeta(t)=0$
and this ensures the continuity of $\zeta$ in $0$. On the other hand, by Lemma~\ref{lemma1} one has that $d$ is strictly decreasing. Thus $\zeta$ is strictly increasing,
and hence an index function.

Because of $0<\frac{1}{q}<1$  this proves the theorem, since $\varphi:=\zeta^{q^*}$, namely $\varphi(0)=0$ and
$\varphi(t)=\left[d \left(\Psi^{-1}(\sigma(t))  \right)\right]^{q^*}$ when $t > 0$,
is an index function, too.}

\begin{remark}
{\rm Easily one can see that the rate function $\left[d \circ \Psi^{-1} \circ \sigma \right]^{q^*}(t)$ in the variational inequality of Theorem~\ref{theo3} tends to zero
 as $t \to 0$ {\sl slower} than the associated rate function $\sigma(t)$ in the variational inequality  of Theorem~\ref{theo2}. Namely,
 taking into account the one-to-one correspondence between large $R>0$ and small $t$ via $\Psi(R)=\sigma(t)$ and
 $\Psi(R)=\frac{d(R)^{q^*}}{R}$ we have for the quotient function
 $$ \frac{\sigma(t)}{\left[d \left(\Psi^{-1}(\sigma(t))  \right)\right]^{q^*}}= \frac{\Psi(R)}{d(R)^{q^*}}=
 \frac{1}{R} \to 0  \quad \mbox{as} \quad R \to \infty,\;\;\mbox{resp.~}\;\; t \to 0\,.$$
As a consequence the situation of approximate source conditions occurring in Theorem~\ref{theo3} leads to lower convergence rates of Tikhonov regularization obtained from Theorem~\ref{theo1} than the situation of canonical source conditions that appears in Theorem~\ref{theo2}.

}\end{remark}

\pagebreak

\begin{example} \label{ex:final}
{\rm Concerning logarithmic rates as an example we are going to conclude the paper with a brief study  that outlines the specific potential of variational inequalities \req{vareqnew}
for extracting both solution smoothness of $\xi$ and nonlinearity conditions on $F$ at $\xdag$ in one index function $\varphi$ which determines the convergence rate.
Let in that example with some $C>0$
\begin{equation} \label{eq:exphi}
\varphi(t)\, = \, \begin{cases}
\;0 & (t=0)\\
\;C\,[\log(1/t)]^{-\mu} & (0<t \le e^{-\mu-1})
\end{cases}
\end{equation}
hold. From that assumption we derive for all $\mu>0$ immediately by Theorem~\ref{theo1} a logarithmic convergence rate
\begin{equation} \label{eq:slowlog}
D_\xi(x_{\alpha(\delta)}^\delta,\xdag)\,=\,{\mathcal
O}\left([\log (1/\delta)]^{-\mu}\right) \quad \mbox{as}\quad \delta \to 0\,,
\end{equation}
which is slower than every power rate \req{prop1rate} for any $\kappa>0$.
Now the function \req{exphi} with slow decay to zero as $t \to 0$ can be a consequence of two completely
different causes characterized by the following two situations (I) and (II), respectively:
\begin{itemize}
\item[(I)] Let $\sigma=\varphi$, i.e., a very weak logarithmic structural condition \req{newstructure} is valid, and let hold the canonical source condition \req{benchmark},
which expresses in our context the strong smoothness assumption on the solution.
Then by Theorem~\ref{theo2} in connection with Theorem~\ref{theo1} we obtain the logarithmic convergence rate \req{slowlog}.

\item[(II)] Let $\sigma(t)=t$, i.e., a structural condition \req{newstructure} is satisfied, which is the strongest in our sense. However, the canonical source condition \req{benchmark}
is strongly violated, which is expressed by a logarithmic decay
$$d(R)= (\log R)^{-\nu} $$
of the corresponding distance function for some $\nu>0$ and all sufficiently large $R>\overline R>0$. However, since we have for all such $R$ and for $\varepsilon>0$ a constant $K>0$ with
$$\Psi(R)= \frac{1}{R(\log R)^{\nu q^*}} \ge \frac{K}{R^{1+\varepsilon}}\,, $$
this implies $\Psi^{-1}(t) \ge \hat K t^{-1/(1+\varepsilon)}$ for some constant $\hat K>0$ and sufficiently small $t>0$. Hence, by Theorem~\ref{theo3} the function $\varphi$ in \req{vareqnew}
attains the form \req{exphi} with $\mu= \nu q ^*.$

\end{itemize}

}\end{example}

\end{document}